\documentclass[12pt,a4paper]{article}
\usepackage{amssymb}
\usepackage{amsmath}
\usepackage{amsfonts}
\usepackage{amstext}
\usepackage[top=1in, bottom=1in,left=0.75in,right=0.75in]{geometry}

\def \qed {\hfill \vrule height6pt width 6pt depth 0pt}

\begin{document}

\title{Unilateral global bifurcation for fourth-order eigenvalue problems with sign-changing weight
\thanks{Research supported by the NSFC (No. 11061030, No. 10971087).}}
\author{{\small  Guowei Dai\thanks{Corresponding author. Tel: +86 931
7971297.\newline
\text{\quad\,\,\, E-mail address}: daiguowei@nwnu.edu.cn (G. Dai).%
}
} \\
%EndAName
{\small Department of Mathematics, Northwest Normal
University, Lanzhou, 730070, PR China}\\
}
\date{}
\maketitle

\begin{abstract}
In this paper,  we shall establish the unilateral global bifurcation result for a class of fourth-order eigenvalue
problems with sign-changing weight.
Under some natural hypotheses on perturbation function,
we show that $\left(\mu_k^\nu,0\right)$ is a bifurcation
point of the above problems and there are two distinct unbounded continua,
$\left(\mathcal{C}_{k}^\nu\right)^+$ and $\left(\mathcal{C}_{k}^\nu\right)^-$,
consisting of the bifurcation branch $\mathcal{C}_{k}^\nu$ from $\left(\mu_k^\nu, 0\right)$,
where $\mu_k^\nu$ is the $k$-th positive or negative eigenvalue of the linear problem corresponding to the above problems,
$\nu\in\{+,-\}$. As the applications of the above result, we study the existence of nodal solutions
for a class of fourth-order eigenvalue problems with sign-changing weight. Moreover, we also establish the Sturm type comparison theorem for fourth-order problems with sign-changing weight.
\end{abstract}\textbf{\ }

\numberwithin{equation}{section}

\numberwithin{equation}{section}

\section{\small Introduction}

\quad\, It is well known that fourth-order elliptic problems arise in many applications, such as Micro Electro
Mechanical systems, thin film theory, surface diffusion on solids, interface dynamics, flow in Hele-Shaw
cells, phase field models of multi-phase systems and the deformation of an elastic beam, see, for example, [\ref{FW}, \ref{M}] and the references therein.
Thus, there are many papers concerning the existence and multiplicity of positive solutions and sign-changing solutions addressed by using
different methods, such as those of topology degree theory, critical point theory, the fixed point theorem in cones and
bifurcation techniques [\ref{BW}, \ref{LO}, \ref{MW}, \ref{MX}, \ref{Y}]. Problems with sign-changing weight arise from the selection-migration model in population genetics. In
this model, weight function $m$ changes sign corresponding to the fact that an allele
$A_1$ holds an advantage over a rival allele $A_2$ at
same points and is at a disadvantage at others; the parameter $\theta$
corresponds to the reciprocal of diffusion, for details, see [\ref{F}].

Recently, Ma et al. [\ref{MG}] established the existence of the principal eigenvalues of the following
linear indefinite weight problem
\begin{equation}
\left\{
\begin{array}{l}
u''''=\lambda g(t)u,\, \ \  \quad\quad\quad  t\in (0,1),\\
u(0)=u(1)=u''(0)=u''(1)=0,
\end{array}
\right.\nonumber
\end{equation}
where $g:[0,1]\rightarrow \mathbb{R}$ is a continuous sign-changing function. They also proved the existence of positive solutions for the corresponding
nonlinear indefinite weight problem. However, there is no any information on the high eigenvalues and the existence of sign-changing solutions for the corresponding
nonlinear indefinite weight problem.

In [\ref{DM1}], Dai and Ma established a Dancer-type unilateral global bifurcation result for one-dimensional $p$-Laplacian problem.
Later, Dai and Ma [\ref{DM2}] established the spectrum of the following eigenvalue problem
\begin{equation}\label{oepp}
\left\{
\begin{array}{l}
u''''=\mu m(t)u, \,\, t\in (0,1),\\
u(0)=u(1)=u''(0)=u''(1)=0,
\end{array}
\right.
\end{equation}
where $\mu$ is a real parameter and $m$ is sign-changing weight.
They proved there exists a unique sequence of eigenvalues for the above problem.
Each eigenvalue is simple, the $k$-th eigenfunction, corresponding to the $k$-th
positive or negative eigenvalue, has exactly $k-1$ generalized simple zeros in $(0,1)$.

In this paper, based the spectral theory of [\ref{DM2}],
we shall establish the similar results to Dai and Ma [\ref{DM1}] about the continua of
solutions for the following fourth-order eigenvalue problem
\begin{equation}\label{oeg}
\left\{
\begin{array}{l}
u ''''=\mu m(t)u+g(t,u,\mu), \ \  t\in (0,1),\\
u(0)=u(1)=u''(0)=u''(1)=0,
\end{array}
\right.
\end{equation}
where $m$ is a sign-changing
function, $g:(0,1)\times \mathbb{R}^2\rightarrow\mathbb{R}$ satisfies the Carath\'{e}odory condition. Let $I:=(0,1)$ and
\begin{equation}
M(I):=\left\{m\in C(\overline{I})\big|\text{meas}\{t\in I,m(t)>0\}\neq0\right\}.\nonumber
\end{equation}
We also assume that the perturbation function $g:I\times \mathbb{R}^2\rightarrow\mathbb{R}$
is continuous and satisfies the following hypotheses:
\begin{equation}\label{ogc}
\lim_{s\rightarrow0}\frac{g(t,s,\mu)}{\vert s\vert}=0
\end{equation}
uniformly for $t\in I$ and $\mu$ on bounded sets.

\indent Under the condition of $m\in M(I)$ and (\ref{ogc}), we shall show that
$\left(\mu_k^\nu,0\right)$ is a bifurcation point of (\ref{oeg}) and there are two distinct unbounded continua,
$\left(\mathcal{C}_{k}^\nu\right)^+$ and $\left(\mathcal{C}_{k}^\nu\right)^-$,
consisting of the bifurcation branch $\mathcal{C}_{k}^\nu$ from $\left(\mu_k^\nu, 0\right)$,
where $\mu_k^\nu$ is the $k$-th positive or negative eigenvalue of
the linear problem corresponding to (\ref{oeg}), where $\nu\in\{+,-\}$.

Based on the above result, we
investigate the existence of nodal solutions for the following fourth-order problem
\begin{equation}\label{oef}
\left\{
\begin{array}{l}
u''''-\gamma m(t)f(u)=0,\,\, t\in I,\\
u(0)=u(1)=u''(0)=u''(1)=0,
\end{array}
\right.
\end{equation}
where $f\in C(\mathbb{R})$, $\gamma$ is a parameter.

The rest of this paper is arranged as follows. In Section 2, we establish the unilateral global bifurcation theory for (\ref{oeg}).
In Section 3, we establish the Sturm type comparison theorem for fourth-order problems with sign-changing weight.
In Section 4, we prove the existence of nodal solutions for (\ref{oef}) under the linear growth condition on $f$.

\section{Unilateral global bifurcation results}
\quad\, We start by considering the following auxiliary problem
\begin{equation}\label{oee}
\left\{
\begin{array}{l}
-u''=e(t),\,\, t\in I,\\
u(0)=u(1)=0
\end{array}
\right.
\end{equation}
for a given $e\in C(\overline{I})$.
It is well known that for every given $e\in C(\overline{I})$ there is a unique solution
$u\in C^2(\overline{I})$ to the problem (\ref{oee}) (see [\ref{BS}]). Let $\Lambda(e)$ denote the unique solution to (\ref{oee})
for a given $e\in C(\overline{I})$. By the results of [\ref{BS}], we can easily show that $\Lambda:C^k(\overline{I})\rightarrow C^{k+2}(\overline{I})$ is continuous for any $k\geq 0$, $k\in \mathbb{N}\cup\{0\}$. Hence, $\Lambda^2: C(\overline{I})\rightarrow C^{3}(\overline{I})$ is compact.

Now, the problem (\ref{oeg}) can be restated as an operator equation
\begin{equation}
u=\mu\Lambda^2(mu).\nonumber
\end{equation}
Define $T_{\mu}(u)=\mu\Lambda^2(mu)$. Let $E = \left\{u\in C^3(\overline{I})\big|u(0)=u(1)=u''(0)=u''(1)=0\right\}$ with the norm
$$\Vert u\Vert=\max_{t\in\overline{I}}\vert u(t)\vert+\max_{t\in\overline{I}}\vert u'(t)\vert+\max_{t\in\overline{I}}\vert u''(t)\vert+\max_{t\in\overline{I}}\vert u'''(t)\vert$$
and $\Psi_{\mu}$ be defined on $E$ by
\begin{equation}
\Psi_{\mu}(u)=u-T_{\mu}(u).\nonumber
\end{equation}
It is no difficult to show that $\Psi_{\mu}$ is a nonlinear
compact perturbation of the identity. Thus the Leray-Schauder degree
$\deg\left(\Psi_{\mu}, B_r(0),0\right)$ is well-defined for
arbitrary $r$-ball $B_r(0)$ and $\mu\neq \mu_k^\nu$.
\\

Firstly, we
can compute $\deg\left(\Psi_{\mu}, B_r(0),0\right)$ for any $r>0$ as follows.
\\ \\
\noindent\textbf{Lemma 2.1.} \emph{For $r>0$, we have
$$\deg\left(\Psi_{\mu}, B_r(0),0\right)
=\left\{\aligned 1, \ \ \ \qquad \ \ & \text{if}\ \mu\in
\left(\mu_1^-,
\mu_1^+\right),\\
(-1)^k, \ \ \ \ \ \ & \text{if}\ \mu\in \left(\mu_k^+, \mu_{k+1}^+\right),\ k\in \mathbb{N},\\
(-1)^k, \ \ \ \ \ \ & \text{if}\ \mu\in \left(\mu_{k+1}^-, \mu_{k}^-\right),\ k\in \mathbb{N}.\\
\endaligned
\right.
$$
}
\noindent\textbf{Proof.} We divide the proof into two cases.

{\it Case 1.} $\mu\geq 0$.

Since $T_{\mu}$ is compact and linear, by [\ref{De}, Theorem 8.10],
\begin{equation}
\deg\left(\Psi_{\mu}, B_r(0),0\right)=(-1)^{m(\mu)},\nonumber
\end{equation}
where $m(\mu)$ is the sum of algebraic multiplicity of the
eigenvalues $\mu$ of (\ref{oepp}) satisfying $\mu^{-1}\mu_k^+<1$.
If $\mu\in [0, \mu^+_1)$, then there are no such $\mu$ at
all, then
\begin{equation}
\deg\left(\Psi_{\mu}, B_r(0),0\right)=(-1)^{m(\mu)}=(-1)^0=1.\nonumber
\end{equation}
If $\mu\in \left(\mu_k^+, \mu_{k+1}^+\right)$ for some $k\in
\mathbb{N}$, then
\begin{equation}
\left(\mu^+_j\right)^{-1}\mu>1, \ \ j\in \{1, \cdots, k\}.\nonumber
\end{equation}
This implies
\begin{equation}
\deg\left(\Psi_{\mu}, B_r(0),0\right)=(-1)^k.\nonumber
\end{equation}

{\it Case 2.} $\mu<0$.

In this case, we consider a new sign-changing eigenvalue problem
$$
\left\{
\begin{array}{l}
u''''-\hat\mu \hat m(t)u=0,\,\,t\in I,\\
u(0)=u(1)=u''(0)=u''(1)=0,
\end{array}
\right.
$$
where $\hat\mu=-\mu$, $\hat m(t)=-m(t)$. It is easy to check
that
\begin{equation}
\hat\mu_k^+=-\mu_k^-, \ \  k\in
\mathbb{N}.\nonumber
\end{equation}
Thus, we may use the result obtained in \emph{Case 1} to deduce the desired
result.  \qed\\

\indent Define the operator $R:\mathbb{R}\times E\rightarrow E$
by
$$
R(\mu,u)(t):=\mu\Lambda^2(mu)+\Lambda^2 g(t,u,\mu).\nonumber
$$
Then it is clear that problem (\ref{oeg}) can be equivalently written as
$$
u=R(\mu,u).
$$
Clearly, $R$ is completely continuous from $\mathbb{R}\times E\rightarrow E$ and
$R(\mu,0)=0$, $\forall \mu\in \mathbb{R}$.
\\ \\
\textbf{Theorem 2.1.} \emph{Assume (\ref{ogc}) holds and $m\in M(I)$. Then $\left(\mu_k^\nu,0\right)$ is a bifurcation
point of (\ref{oeg}) and the associated bifurcation branch $\mathcal{C}_k^\nu$ in $\mathbb{R}\times E$
whose closure contains ($\mu_k^\nu, 0$) is either unbounded or contains a pair ($\overline{\mu}, 0$)
with $\overline{\mu}$ is an eigenvalue of (\ref{oepp}) and $\overline{\mu}\neq\mu_k^\nu$.}
\\ \\
\textbf{Proof.} We only prove the case of $\mu_k^+$ because the case of $\mu_k^-$
is similar. From now on, for simplicity, we write $\mu_k = \mu_k^+$. Suppose that $(\mu_k, 0)$
is not a bifurcation point of problem (\ref{oeg}). Then
there exist $\varepsilon > 0$, $\rho_0> 0$ such that for $\vert \mu-\mu_k\vert\leq\varepsilon$
and $0<\rho < \rho_0$ there is no
nontrivial solution of the equation
\begin{equation}
u-R(\mu,u)=0\nonumber
\end{equation}
with $\Vert v\Vert=\rho$. From the invariance of the degree under a compact
homo-topology we obtain that
\begin{equation}\label{edc}
\text{deg}\left(I-R(\mu,\cdot),B_\rho(0),0\right)\equiv constant
\end{equation}
for $\mu\in[\mu_k-\varepsilon,\mu_k+\varepsilon]$.

By taking $\varepsilon$ smaller if necessary, we can assume that there is no eigenvalue
of (\ref{oepp}) in $(\mu_k,\mu_k+\varepsilon]$. Fix $\mu\in(\mu_k,\mu_k+\varepsilon]$.
We claim that the equation
\begin{equation}\label{es}
u-\left(\mu\Lambda^2(mu)+s\Lambda^2 g(t,u,\mu)\right)=0
\end{equation}
has no solution $u$ with $\Vert u\Vert=\rho$ for every $s\in\overline{I}$ and $\rho$ sufficiently small.
Suppose on the contrary, let $\{u_n\}$ be the solution of (\ref{es}) with $\Vert u_n\Vert\rightarrow 0$
as $n\rightarrow+\infty$.

Let $v_n:=u_n/\Vert u_n\Vert$,  then $v_n$ should be a solution of problem
\begin{eqnarray}\label{evs0}
v_n(t)&=&\Lambda^2\left(\mu mv_n+s \frac{g(t,u_n,\mu)}{\Vert u_n\Vert}\right).
\end{eqnarray}
Let
\begin{equation}
\widetilde{g}(t,u,\mu)=\max_{0\leq \vert s\vert\leq u}\vert g(t,s,\mu)\vert\,\,
\text{for}\,\, t\in I \text{\,\,and\,\,} \mu \text{\,\,on bounded sets},\nonumber
\end{equation}
then $\widetilde{g}$ is nondecreasing with respect to $u$ and
\begin{equation}\label{eg0+}
\lim_{ u\rightarrow 0^+}\frac{\widetilde{g}(t,u,\mu)}{u}=0
\end{equation}
uniformly for $t\in I$ and $\mu$ on bounded sets.
Further it follows from (\ref{eg0+}) that
\begin{equation}\label{egn0}
\frac{g(t,u,\mu)}{\Vert u\Vert} \leq\frac{
\widetilde{g}(t,\vert u\vert,\mu)}{\Vert u\Vert} \leq \frac{
\widetilde{g}(t,\Vert u\Vert_\infty,\mu)}{\Vert u\Vert} \leq
\frac{ \widetilde{g}(t,\Vert u\Vert,\mu)}{\Vert
u\Vert}\rightarrow0\ \  \text{as}\,\, \Vert
u\Vert\rightarrow 0
\end{equation}
uniformly for $t\in I$ and $\mu$ on bounded sets.

By (\ref{evs0}), (\ref{egn0}) and compactness of $\Lambda^2$, we obtain that for
some convenient subsequence $v_n\rightarrow v_0$ as $n\rightarrow+\infty$. Now $v_0$
verifies the equation
\begin{equation}
v_0''''= \mu mv_0\nonumber
\end{equation}
and $\Vert v_0\Vert = 1$. This implies that $\mu$ is an eigenvalue of (\ref{oepp}).
This is a contradiction.
From the invariance of the degree under
homo-topology and Lemma 2.1 we then obtain
\begin{equation}\label{edFk}
\deg\left(I-R(\mu,\cdot), B_r(0),0\right)=\deg\left(\Psi_{\mu}, B_r(0),0\right)=(-1)^k.
\end{equation}
Similarly, for $\mu\in [\mu_k-\varepsilon, \mu_k)$ we find that
\begin{equation}\label{edFk1}
\deg\left(I-R(\mu,\cdot), B_r(0),0\right)=(-1)^{k-1}.
\end{equation}
Relations (\ref{edFk}) and (\ref{edFk1}) contradicts (\ref{edc}) and hence $(\mu_k, 0)$ is a
bifurcation point of problem (\ref{oeg}).

By standard arguments in global bifurcation theory (see [\ref{R2}]), we can show
the existence of a global branch of solutions of (\ref{oeg}) emanating from
$(\mu_k, 0 )$.\qed\\

Now, we give the definitions of nodal solution, generalized simple zero and generalized double zero.
\\ \\
\textbf{Definition 2.1.} Let $u$ be a nontrivial solution of (\ref{oeg}) and $t_*$ be a zero of $u$.
We call that $t_*$ is a generalized simple zero if $u''(t_*)=0$ but $u'(t_*)\neq 0$ or $u'''(t_*)\neq 0$.
Otherwise, we call that $t_*$ is a generalized double zero. If there is no generalized double zero of $u$, we call that $u$ is a nodal solution.
\\

Next, we prove that the first choice of the alternative of Theorem 2.1 is the only possibility.
Let $S_k^+$ denote the set of
functions in $E$ which have exactly $k-1$ generalized simple zeros in $I$ and are positive near $t=0$, and set $S_k^-=-S_k^+$, and
$S_k =S_k^+\cup S_k^-$. Clearly, they are disjoint and open in $E$. Finally, let
$\Phi_k^{\pm}=\mathbb{R}\times S_k^{\pm}$ and $\Phi_k=\mathbb{R}\times S_k$ under the
product topology.
\\ \\
\textbf{Lemma 2.2.} \emph{If $(\mu, u)$ is a solution of (\ref{oeg})
and $u$ has a generalized double zero, then $u \equiv 0$.}
\\ \\
\textbf{Proof.} Let $u$ be a solution of (\ref{oeg}) and $t^*\in\overline{I}$ be a generalized double zero, i.e., $u(t_*)=u'(t_*)=u''(t_*)=u'''(t_*)=0$.
We note that
\begin{equation}
u(t)=\int_{t_*}^t\int_{t_*}^s\int_{t_*}^\tau\int_{t_*}^\rho\left(\mu m(\xi)u(\xi)+g(\xi,u(\xi),\mu)\right)\,d\xi d\rho d\tau ds.\nonumber
\end{equation}
First, we consider $t\in[0, t^*]$. Then
\begin{eqnarray}
\vert u(t)\vert&=&\left\vert\int_{t_*}^t\int_{t_*}^s\int_{t_*}^\tau\int_{t_*}^\rho\left(\mu m(\xi)u(\xi)+g(\xi,u(\xi),\mu)\right)\,d\xi d\rho d\tau ds\right\vert\nonumber\\
&\leq&\int_{t}^{t_*}\left\vert\int_{t_*}^\tau\int_{t_*}^\rho\left(\mu m(\xi)u(\xi)+g(\xi,u(\xi),\mu)\right)\,d\xi d\rho\right\vert d\tau\nonumber\\
&\leq&\int_{t}^{t_*}\int_{\tau}^{t_*}\left\vert\left(\mu m(\xi)u(\xi)+g(\xi,u(\xi),\mu)\right)\right\vert\,d\xi d\tau\nonumber\\
&\leq&\int_{t}^{t_*}\left\vert\left(\mu m(\xi)u(\xi)+g(\xi,u(\xi),\mu)\right)\right\vert\,d\xi,\nonumber
\end{eqnarray}
furthermore,
\begin{eqnarray}
\vert u(t)\vert&\leq&\int_{t}^{t_*}\left\vert\left(\mu m(\tau)u(\tau)+g(\tau,u(\tau),\mu)\right)\right\vert\,d\tau \nonumber\\
&\leq&\int_t^{t^*}
 \left\vert\mu m(\tau)
+\frac{g(\tau,u(\tau),\mu)}{u(\tau)}\right\vert u(\tau)\,d\tau \nonumber\\
&\leq&\int_t^{t^*}
 \left(\mu \vert m(\tau)\vert
+\left\vert\frac{g(\tau,u(\tau),\mu)}{u(\tau)}\right\vert\right)\vert u(\tau)\vert\,d\tau.\nonumber
\end{eqnarray}
In view of (\ref{ogc}), for any $\varepsilon>0$, there exists a constant $\delta>0$ such that
\begin{equation}
\vert g(t,s,\mu)\vert\leq \varepsilon \vert s\vert\nonumber
\end{equation}
uniformly with respect to $t\in I$ and fixed $\mu$ when $\vert s\vert\in[0,\delta]$.
Hence,
\begin{equation}
\vert u(t)\vert\leq \int_t^{t^*}
\left( \mu\vert m(\tau)\vert+\varepsilon+\max_{s\in\left[\delta,\Vert u\Vert_\infty\right]}
\left\vert\frac{g(\tau,s,\mu)}{s}\right\vert\right) \vert u(\tau)\vert\,d\tau.\nonumber
\end{equation}
By Gronwall-Bellman inequality [\ref{Bre}], we get $u \equiv 0$ on $[0, t^*]$.
Similarly, we also can get $u \equiv 0$ on $[t^*, 1]$
and the proof is complete.\qed
\\ \\
\textbf{Lemma 2.3.} \emph{The last alternative of Theorem 2.1 is
impossible if
$\mathcal{C}_k^\nu\subset\Phi_k\cup\{(\mu_k^\nu,0)\}$.}
\\ \\
\textbf{Proof.} Suppose on the contrary, if there exists $(\mu_m,u_m)\rightarrow\left(\mu_j^\nu,0\right)$
when $m\rightarrow+\infty$ with $(\mu_m,u_m)\in \mathcal{C}_k^\nu$, $u_m \not\equiv 0$ and $j\neq k$.
Let $w_m:=u_m/\Vert u_m\Vert$, then $w_m$ should be a solution of problem
\begin{eqnarray}\label{evgm}
w(t)&=&\Lambda^2\left(\mu mw+\frac{g(t,u_m,\mu)}{\Vert u_m\Vert}\right).
\end{eqnarray}
By (\ref{egn0}), (\ref{evgm}) and the compactness of $\Lambda^2$ we obtain that
for some convenient subsequence $w_m\rightarrow w_0$ as $m\rightarrow+\infty$. Now $w_0$ verifies the equation
\begin{equation}
w_0''''= \mu_j^\nu m(t)w_0\nonumber
\end{equation}
and $\Vert w_0\Vert = 1$. Hence $w_0\in S_j$ which is an open set in $E$, and as a
consequence for some $m$ large enough, $u_m \in S_j$, and this is a contradiction.\qed
\\ \\
\textbf{Theorem 2.2.} \emph{Assume (\ref{ogc}) holds and $m\in M(I)$, then from each $\left(\mu_k^\nu,0\right)$ it
bifurcates an unbounded continuum $\mathcal{C}_k^\nu$ of solutions to problem (\ref{oeg}),
with exactly $k-1$ simple zeros, where $\mu_k^\nu$ is the eigenvalue of problem (\ref{oepp}).}\\ \\
\textbf{Proof.}
Taking into account Theorem 2.1 and Lemma 2.3, we only need to prove that
$\mathcal{C}_k^\nu\subset\Phi_k\cup\{(\mu_k^\nu,0)\}$.
Suppose $\mathcal{C}_k^\nu\not\subset\Phi_k\cup\left\{\left(\mu_k^\nu,0\right)\right\}$.
Then there exists
$(\mu, u)\in \mathcal{C}_k^\nu\cap(\mathbb{R}\times \partial S_k)$ such that $(\mu, u) \neq
\left(\mu_k^\nu, 0\right)$ and $(\mu_n,u_n)\rightarrow(\mu, u)$ with $(\mu_n,u_n)
\in \mathcal{C}_k^\nu \cap(\mathbb{R}\times S_k)$.
Since $u\in \partial S_k$, by Lemma 2.2, $u\equiv 0$. Let $v_n:=u_n/\Vert u_n\Vert$, then $v_n$ should be a solution of problem
\begin{eqnarray}\label{ewgm}
v(t)&=&\Lambda^2\left(\mu mv+\frac{g(t,u_n,\mu)}{\Vert u_n\Vert}\right)
\end{eqnarray}
By (\ref{egn0}), (\ref{ewgm}) and the compactness of $\Lambda^2$ we obtain that for some
convenient subsequence $v_n\rightarrow v_0$ as $n\rightarrow+\infty$. Now $v_0$ verifies the equation
\begin{equation}
v_0''''= \mu m(t)v_0\nonumber
\end{equation}
and $\Vert v_0\Vert = 1$. Hence $\mu = \mu_j^\nu$, for some $j \neq k$.
Therefore, $(\mu_n,u_n)\rightarrow\left(\mu_j^\nu, 0\right)$ with $(\mu_n,u_n)\in \mathcal{C}_k^\nu\cap (\mathbb{R}\times S_k)$.
This contradicts Lemma 2.3.\qed\\

Using the similar method to prove [\ref{DM1}, Theorem 3.2] with obvious changes, we may obtain the
following result.\\ \\
\textbf{Theorem 2.3.} \emph{Assume (\ref{ogc}) holds and $m\in M(I)$, then there are two distinct unbounded continua,
$\left(\mathcal{C}_{k}^\nu\right)^+$ and $\left(\mathcal{C}_{k}^\nu\right)^-$,
consisting of the bifurcation branch $\mathcal{C}_{k}^\nu$. Moreover, for $\sigma\in\{+,-\}$, we have}
$$
\left(\mathcal{C}_{k}^\nu\right)^\sigma\subset \left(\{\left(\mu_k^\nu,0\right)\}\cup\left(\mathbb{R}\times S_k^\sigma\right)\right).
$$\nonumber
%以上结果可以推到p-双调和，但下面的比较定理难推广
\section{Sturm type comparison theorem}

\quad\, In this section, we shall establish the Sturm type comparison theorem for fourth-order differential equations with sign-changing weight, which will be used later.\\ \\
\textbf{Lemma 3.1.} \emph{Let $b_2(t)>b_1(t)>0$ for $t\in I$ and $b_i(t)\in C(\overline{I})$, $i=1,2$. Also let $u_1$, $u_2$
be solutions of the following differential equations:
\begin{equation}
u''''=b_i(t)u, \,\, t\in I, i=1,2,\nonumber
\end{equation}
respectively. If $u_1$ has $k$ generalized simple zeros in $I$, then $u_2$ has at least $k+1$ generalized simple zeros in $I$.}
\\ \\
\textbf{Proof.} Let $c$ and $d$ be any two consecutive generalized simple zeros of $u_1$ in $\overline{I}$. Then we can assume without loss of generality that
$u_1(t)>0$, $u_2(t)>0$ in $(c,d)$. Then an easy calculation shows that
$$
\int_c^d\left(u_1''''u_2-u_2''''u_1\right)\,dt=\int_c^d\left(b_1-b_2\right)u_1u_2\,dt<0.\eqno (3.1)
$$
The left-hand side of (3.1) equals
\begin{equation}
u_1'''(d)u_2(d)-u_1'''(c)u_2(c)+u_1'(d)u_2''(d)-u_1'(c)u_2''(c).\nonumber
\end{equation}
Next, we shall show that
\begin{equation}
u_1'''(d)u_2(d)-u_1'''(c)u_2(c)+u_1'(d)u_2''(d)-u_1'(c)u_2''(c)\geq 0.\nonumber
\end{equation}
In fact, if this occurs, we arrive a contradiction.
We divide the proof into two steps.

\emph{Step 1}: We show that $u_1'''(d)u_2(d)-u_1'''(c)u_2(c)\geq 0$.
Let $v:=u_1''$. We consider the system:
$$
\left\{
\begin{array}{l}
u_1''=v, \,\, t\in I,\\
v''=b_1u_1.
\end{array}
\right.
$$
By simple computation, one has
$$
u_1'v'=\frac{v^2}{2}+\frac{bu_1^2}{2}+C\eqno (3.2)
$$
for any constant $C$.
Let $t_0\in(c,d)$ be the point satisfying
\begin{equation}
u_1(t_0)=\max_{t\in[c,d]}u_1(t).\nonumber
\end{equation}
Then (3.2) implies
$$
0=\frac{v^2(t_0)}{2}+\frac{bu_1^2(t_0)}{2}+C.
$$
It follows $C<0$. Putting $c$ into (3.2), we have
$$
u_1'(c)u_1'''(c)=C<0.
$$
Using this and the fact $u_1'(c)\geq 0$, we get $u_1'''(c)<0$. Similarly, we can show that $u_1'''(d)>0$.
Hence, we have $u_1'''(d)u_2(d)-u_1'''(c)u_2(c)\geq 0$.

\emph{Step 2}: We show that $u_1'(d)u_2''(d)-u_1'(c)u_2''(c)\geq 0$.

It suffices to show that $u_2''(c)\leq 0$ and $u_2''(d)\leq 0$ since the facts $u_1'(c)\geq 0$ and $u_1'(d)\leq0$.
Suppose on the contrary that $u_2''(c)> 0$ or $u_2''(d)> 0$, we shall deduce a contradiction.

Let $u_*:=u_2(t)+1$. Then $u_*''''=b_2u_2$ and $u_*\geq1$ in $(c,d)$. For some $\varepsilon>0$ small enough, let $\widetilde{u}\in C^4([-\varepsilon,1+\varepsilon])$ and
$\widetilde{b}\geq 0$ be such that
$\widetilde{u}(-\varepsilon)=\widetilde{u}(1+\epsilon)=\widetilde{u}''(-\varepsilon)=\widetilde{u}''(1+\varepsilon)=0$,
$\widetilde{u}\big|_I=u_*$ and $\widetilde{u}''''=\widetilde{b}\widetilde{u}$. Then we have
$$
\left\{
\begin{array}{l}
\widetilde{u}''''=\widetilde{b} \widetilde{u}, \,\, t\in (-\varepsilon,1+\varepsilon),\\
\widetilde{u}(-\varepsilon)=\widetilde{u}(1+\epsilon)=\widetilde{u}''(-\varepsilon)=\widetilde{u}''(1+\varepsilon)=0.
\end{array}
\right.
$$
Set $a:=(c-\varepsilon)/2$ and $b:=(d+1+\varepsilon)/2$.
Let $\overline{u}\in C^4([a,b])$ and $\overline{b}\geq0$ be such that $\overline{u}\big|_I=\widetilde{u}$, $\overline{u}\geq0$ in ($a,b$) and
$\overline{u}(a)=\overline{u}(b)=\overline{u}''(a)=\overline{u}''(b)=0$ and
$\overline{u}''''=\overline{b}\overline{u}$. Set $w:=\overline{u}''$, then $w$ should be a solution of the problem
$$
\left\{
\begin{array}{l}
w''=\overline{b}\overline{u}, \,\, t\in (a,b),\\
w(a)=w(b)=0.
\end{array}
\right.
$$
The Strong Maximum Principle implies that $w<0$ in $(a,b)$. This follows that $u_2''\leq 0$ in $[c,d]$.\qed \\

Let
\begin{equation}
I^+:=\left\{t\in \overline{I}\,|\, m(t)>0\right\}, \ \ I^-:=\left\{t\in \overline{I}\,|\, m(t)<0\right\}.\nonumber
\end{equation}
\noindent\textbf{Lemma 3.2.} \emph{Assume $m\in M(I)$. Let $\widehat{I}=(a,b)$ be such
that $\widehat{I}\subset I^+$ and
\begin{equation}
\text{meas}\,\widehat{I}>0.
\nonumber
\end{equation}
Let $g_n:I\to (0, +\infty)$ be continuous function and such that
\begin{equation}
\lim_{n\to +\infty} g_n(t)=+\infty\ \ \text{uniformly on}\
\widehat{I}.\nonumber
\end{equation}
Let $y_n$ be a solution of the equation
\begin{equation}
\left\{
\begin{array}{l}
y_n''''=m(t)g_n(t)y_n, \,\, t\in I,\\
u(0)=u(1)=u''(1)=u''(1)=0.
\end{array}
\right.\nonumber
\end{equation}
Then the number of zeros of $y_n$ in $I$ goes to infinity as $n\to +\infty$.}
\\ \\
\textbf{Proof.} After taking a subsequence if
necessary, we may assume that
\begin{equation}
m(t)g_{n_j}(t)\geq \lambda_j, \ \  t\in \widehat{I}\nonumber
\end{equation}
as $j\to +\infty$,
where $\lambda_j$ is the $j$-th eigenvalue of the following problem
\begin{equation}
u''''=\lambda u(t),\,\,t\in I.\nonumber
\end{equation}
Let $\varphi_j$ be the corresponding eigenvalue of $\lambda_j$.
It is easy to check that the distance between any
two consecutive zeros of $\varphi_j$ is $1/j$ (also see [\ref{DO}]).
Hence, the number of zeros of $\varphi_j\big|_{\widehat{I}}$ goes to infinity as $j\to +\infty$. By Lemma 3.1,
one has that the number of zeros of $y_n|_{\widehat{I}}$ goes to infinity as $n\to +\infty$.
It follows the desired results.\qed\\

Similarly, we also have:\\
\\
\noindent\textbf{Lemma 3.3.} \emph{Assume $m\in M(I)$. Let $\widetilde{I}=(c,d)$ be such
that $\widetilde{I}\subset I^-$ and
\begin{equation}
\text{meas}\,\widetilde{I}>0.
\nonumber
\end{equation}
Let $g_n:I\to (-\infty, 0)$ be continuous function and such that
\begin{equation}
\lim_{n\to +\infty} g_n(t)=-\infty \ \ \text{uniformly on}\
\widetilde{I}.\nonumber
\end{equation}
Let $y_n$ be a solution of the equation
\begin{equation}
y_n''''=m(t)g_n(t)y_n, \,\, t\in I.\nonumber
\end{equation}
Then the number of zeros of $y_n$ goes to infinity as $n\to +\infty$.}

\section{Existence of nodal solutions of (\ref{oef})}

\quad\, In this section, we shall investigate the existence and multiplicity of nodal
solutions to the problem (\ref{oef}) under the linear growth condition on $f$.\\

\indent Firstly, we suppose that\\

\emph{(${H}_1$) $f\in C(\mathbb{R},\mathbb{R})$ with
$f(s)s>0$ for $s\neq0$;}

\emph{($H_2$) there exist $f_0$, $f_\infty\in (0, +\infty)$ such that}
\begin{equation}
f_0=\lim_{\vert s\vert\rightarrow0}\frac{f(s)}{s},\,\, f_\infty=\lim_{\vert
s\vert\rightarrow+\infty}\frac{f(s)}{s}.\nonumber
\end{equation}

\indent Let $\mu_k^\pm$ be the $k$-th positive or negative eigenvalue of (\ref{oepp}).
Applying Theorem 2.3, we shall establish the existence of nodal solutions of (\ref{oef}) follows.\\ \\
\textbf{Theorem 4.1.} \emph{Let ($H_1$), ($H_2$) hold and $m\in M(I)$.
Assume that for some $k \in \mathbb{N}$, either
\begin{equation}
\gamma\in\left(\frac{\mu_k^+}{f_\infty},\frac{\mu_k^+}{f_0}\right)
\cup\left(\frac{\mu_k^-}{f_0},\frac{\mu_k^-}{f_\infty}\right)\nonumber
\end{equation}
or
\begin{equation}
\gamma\in\left(\frac{\mu_k^+}{f_0},\frac{\mu_k^+}{f_\infty}\right)
\cup\left(\frac{\mu_k^-}{f_\infty},\frac{\mu_k^-}{f_0}\right).\nonumber
\end{equation}
Then (\ref{oef}) has two solutions $u_k^+$ and $u_k^-$ such
that $u_k^+$ has exactly $k-1$ generalized simple zeros in $I$ and is positive near 0,
and $u_k^-$ has exactly $k-1$ generalized simple zeros in $I$ and is negative near 0.}
\\ \\
\textbf{Proof.} We only prove the case of $\gamma>0$. The case of
$\gamma<0$ is similar. Consider the problem
$$
\left\{
\begin{array}{l}
u''''=\mu\gamma m(t)f(u),\,\, t\in I,\\
u(0)=u(1)=u''(0)=u''(1)=0.
\end{array}
\right.\eqno (4.1)
$$
Let $\zeta\in C(\mathbb{R})$ be such that
\begin{equation}
f(u)=f_0u+\zeta(u)\nonumber
\end{equation}
with
\begin{equation}\label{eac}
\lim_{\vert u\vert\rightarrow0}\frac{\zeta(u)}{u}=0.\nonumber
\end{equation}
Hence, the condition
(\ref{ogc}) holds. Using Theorem 2.3, we have that
there are two distinct unbounded continua, $\left(\mathcal{C}_k^\nu\right)^+$ and $\left(\mathcal{C}_k^\nu\right)^-$,
consisting of the bifurcation branch $\mathcal{C}_k^\nu$ from $\left(\mu_k^\nu/\gamma f_0, 0\right)$
, such that
\begin{equation}
\left(\mathcal{C}_{k}^\nu\right)^\sigma\subset \left(\left\{\left(\mu_k^\nu,0\right)\right\}\cup\left(\mathbb{R}\times S_k^\sigma\right)\right).\nonumber
\end{equation}

It is clear that any solution of (4.1) of the form $(1, u)$ yields a
solutions $u$ of (\ref{oef}). We shall show that $\left(C_k^+\right)^\sigma$ crosses
the hyperplane $\{1\}\times E$ in $\mathbb{R}\times E$. To  this
end, it will be enough to show that $\left(C_k^+\right)^\sigma$ joins
$\left(\frac{\mu_k^+}{\gamma f_0}, 0\right)$ to
$\left(\frac{\mu_k^+}{\gamma f_\infty}, +\infty\right)$. Let
$\left(\mu_n, y_n\right) \in \left(C_k^+\right)^\sigma$ satisfy $\mu_n+\Vert y_n\Vert\rightarrow+\infty.$
We note that $\mu_n >0$ for all $n \in \mathbb{N}$ since $(0, 0)$ is
the only solution of (4.1) for $\mu = 0$ and
$\left(C_k^+\right)^\sigma\cap\left(\{0\}\times E\right)=\emptyset$.

\emph{Case 1}: $\mu_k^+/f_\infty<\gamma<\mu_k^+/f_0$.

In this case, we only need to show that
\begin{equation}
\left(\frac{\mu_k^+}{\gamma f_\infty},\frac{\mu_k^+}{\gamma
f_0}\right) \subseteq\left\{\mu\in\mathbb{R}:(\mu,u)\in \left(C_k^+\right)^\sigma\right\}.\nonumber
\end{equation}
We divide the proof into two steps.

\emph{Step 1}: We show that if there exists a constant $M>0$ such that
\begin{equation}
\mu_n\subset(0,M]\nonumber
\end{equation}
for $n\in \mathbb{N}$ large enough, then $(C_k^+)^\sigma$ joins
$\left(\mu_k^+/\gamma f_0, 0\right)$ to
$\left(\mu_k^+/\gamma f_\infty, +\infty\right)$.

In this case it follows that
\begin{equation}
\Vert y_n\Vert\rightarrow+\infty.\nonumber
\end{equation}
Let $\xi\in C(\mathbb{R})$ be such that
\begin{equation}
f(u)=f_\infty u+\xi(u).\nonumber
\end{equation}
Then
\begin{equation}
\lim_{\vert u\vert\rightarrow+\infty}\frac{\xi(u)}{u}=0.\nonumber
\end{equation}
Let
\begin{equation}
\widetilde{\xi}(u)=\max_{0\leq \vert s\vert\leq u}\vert \xi(s)\vert.\nonumber
\end{equation}
Then $\widetilde{\xi}$ is nondecreasing and
$$
\lim_{u\rightarrow +\infty}\frac{\widetilde{\xi}(u)}{\vert u\vert}=0.\eqno(4.2)
$$

We divide the equation
\begin{equation}
y_n''''=\mu_n\gamma m(t)f_\infty y_n+\mu_n\gamma m(t)\xi(y_n)\nonumber
\end{equation}
by $\Vert y_n\Vert$ and set $\overline{y}_n = y_n/\Vert y_n\Vert$. Since $\overline{y}_n$ is bounded in $E$,
after taking a subsequence if
necessary, we have that $\overline{y}_n \rightharpoonup \overline{y}$ for some $\overline{y} \in E$. Moreover, from
(4.2) and the fact that $\widetilde{\xi}$ is nondecreasing, we have that
\begin{equation}
\lim_{n\rightarrow+\infty}\frac{ \xi(y_n(t))}{\Vert y_n\Vert}=0\nonumber
\end{equation}
since
\begin{equation}
\frac{ \xi(y_n(t))}{\Vert y_n\Vert}\leq\frac{ \widetilde{\xi}(\vert y_n(t)\vert)}{\Vert y_n\Vert}
\leq\frac{ \widetilde{\xi}(\Vert y_n(t)\Vert_\infty)}{\Vert y_n\Vert}
\leq\frac{ \widetilde{\xi}(\Vert y_n(t)\Vert)}{\Vert y_n\Vert}.\nonumber
\end{equation}
By the continuity and compactness of $\Lambda^2$, it follows that
\begin{equation}
\overline{y}''''=\overline{\mu}\gamma m(t)f_\infty\overline{y},\nonumber
\end{equation}
where
$\overline{\mu}=\underset{n\rightarrow+\infty}\lim\mu_n$, again
choosing a subsequence and relabeling if necessary.

We claim that
\begin{equation}
\overline{y}\in \left(C_k^+\right)^\sigma.\nonumber
\end{equation}

It is clear that $\overline{y}\in \overline{\left(\mathcal{C}_{k}^+\right)^\sigma}\subseteq
\left(\mathcal{C}_{k}^+\right)^\sigma$ since $\left(\mathcal{C}_{k}^+\right)^\sigma$ is closed in $\mathbb{R}\times E$.
Hence, $\overline{\mu}\gamma f_\infty=\mu_k^+$, so
that
\begin{equation}
\overline{\mu}=\frac{\mu_k^+}{\gamma f_\infty}.\nonumber
\end{equation}
Therefore, $(C_k^+)^\sigma$ joins $\left(\mu_k^+/\gamma
f_0, 0\right)$ to $\left(\mu_k^+/\gamma f_\infty,
+\infty\right)$.

\emph{Step 2}: We show that there exists a constant $M$ such that $\mu_n
\in(0,M]$ for $n\in \mathbb{N}$ large enough.

On the contrary, we suppose that
\begin{equation}
\lim_{n\rightarrow +\infty}\mu_n=+\infty.\nonumber
\end{equation}
Since  $\left(\mu_n, y_n\right) \in (C_k^+)^\sigma$, it follows that
\begin{equation}
y_n''''=\gamma\mu_n
m(t)\widetilde{f}_n(t)\varphi(y_n),\nonumber
\end{equation}
where
\begin{equation}
\widetilde{f}_n(t)=\left\{
\begin{array}{l}
\frac{f(y_n(t))}{y_n(t)},\,\, \text{if}\,\,y_n(t)\neq0,\\
f_0,\,\,\quad\,\,\,\,\,\,\text{if}\,\,y_n(t)=0.
\end{array}
\right.\nonumber
\end{equation}
Conditions ($H_1$) and ($H_2$) imply that there exists a positive constant $\varrho$ such that $\widetilde{f}_n(t)>\varrho$
for any $t\in\overline{I}$ and all $n\in \mathbb{N}$.
Then Lemma 3.2 follows that $y_n$ has more than $k$ zeros in $I$ for $n$
large enough, and this contradicts the fact that $y_n$ has exactly $k-1$ zeros in $I$.

\emph{Case 2}: $\mu_k^+/f_0<\gamma<\mu_k^+/f_\infty$.

In this case, we have that
\begin{equation}
\frac{\mu_k^+}{\gamma f_0}<1<\frac{\mu_k^+}{\gamma
f_\infty}.\nonumber
\end{equation}
Assume that $\left(\mu_n, y_n\right) \in \left(C_k^+\right)^\sigma$ is such that
\begin{equation}
\lim_{n\rightarrow+\infty}\left(\mu_n+\Vert y_n\Vert\right)=+\infty.\nonumber
\end{equation}
In view of \emph{Step 2} of \emph{Case 1}, we have known that there exists $M>0$, such that
for $n \in \mathbb{N}$ sufficiently large,
\begin{equation}
\mu_n\in (0,M].\nonumber
\end{equation}
Applying the same method used in \emph{Step 1} of \emph{Case 1}, after
taking a subsequence and relabeling if necessary, it follows that
\begin{equation}
(\mu_n,y_n)\rightarrow\left(\frac{\mu_k^+}{\gamma
f_\infty},+\infty\right) \,\, \text{as}\,\, n\rightarrow+\infty.\nonumber
\end{equation}
Thus, $(C_k^+)^\sigma$ joins $\left(\mu_k^+/\gamma f_0,0\right)$ to $\left(\mu_k^+/\gamma f_\infty,+\infty\right)$.\qed
\\

\indent Using the similar proof with the proof Theorem 4.1, we can obtain the more general results as follows.
\\ \\
\textbf{Theorem 4.2.} \emph{Let ($H_1$), ($H_2$) hold and $m\in M(I)$.
Assume that for some $k, n \in \mathbb{N}$ with $k\leq n$, either
\begin{equation*}
\gamma\in\left(\frac{\mu_n^+}{f_\infty},\frac{\mu_k^+}{f_0}\right)
\cup\left(\frac{\mu_k^-}{f_0},\frac{\mu_n^-}{f_\infty}\right)
\end{equation*}%
or
\begin{equation}
\gamma\in\left(\frac{\mu_n^+}{f_0},\frac{\mu_k^+}{f_\infty}\right)
\cup\left(\frac{\mu_k^-}{f_\infty},\frac{\mu_n^-}{f_0}\right).\nonumber
\end{equation}
Then (\ref{oef}) has $n-k+1$ pairs solutions $u_j^+$ and $u_j^-$ for $j\in\{k,\cdots,n\}$ such
that $u_j^+$ has exactly $j-1$ generalized simple zeros in $I$ and is positive near 0,
and $u_j^-$ has exactly $j-1$ generalized simple zeros in $I$ and is negative near 0.}
\\ \\
\textbf{Remark 4.1.} Clearly, Theorem 1.1 of [\ref{MG}] is the corollary of Theorem 4.1.

\end{document}